\documentclass[10pt]{article}
\usepackage{amsfonts,amsmath,amssymb,cite}
\newcommand{\Real}{\mathop{\textrm{Re}}\nolimits}

\begin{document}
\title{Some integrals and series involving the Gegenbauer polynomials
and the Legendre functions on the cut $(-1,1)$}
\author{Rados{\l}aw Szmytkowski \\*[3ex]
Atomic Physics Division,
Department of Atomic Physics and Luminescence, \\
Faculty of Applied Physics and Mathematics,
Gda{\'n}sk University of Technology, \\
Narutowicza 11/12, PL 80--233 Gda{\'n}sk, Poland \\
email: radek@mif.pg.gda.pl}
\date{\today}
\maketitle
\begin{abstract}
We use the recent findings of Cohl [arXiv:1105.2735] and evaluate two
integrals involving the Gegenbauer polynomials:
$\int_{-1}^{x}\mathrm{d}t\:
(1-t^{2})^{\lambda-1/2}(x-t)^{-\kappa-1/2}C_{n}^{\lambda}(t)$ and
\linebreak
$\int_{x}^{1}\mathrm{d}t\:(1-t^{2})^{\lambda-1/2}(t-x)^{-\kappa-1/2}
C_{n}^{\lambda}(t)$, both with $\Real\lambda>-\frac{1}{2}$,
$\Real\kappa<\frac{1}{2}$, \mbox{$-1<x<1$}. The results are expressed
in terms of the on-the-cut associated Legendre functions
$P_{n+\lambda-1/2}^{\kappa-\lambda}(\pm x)$ and
$Q_{n+\lambda-1/2}^{\kappa-\lambda}(x)$. In addition, we find
closed-form representations of the series
\mbox{$\sum_{n=0}^{\infty}(\pm)^{n}[(n+\lambda)/\lambda]
P_{n+\lambda-1/2}^{\kappa-\lambda}(\pm x)C_{n}^{\lambda}(t)$} and 
\mbox{$\sum_{n=0}^{\infty}(\pm)^{n}[(n+\lambda)/\lambda]
Q_{n+\lambda-1/2}^{\kappa-\lambda}(\pm x)C_{n}^{\lambda}(t)$}, both 
with $\Real\lambda>-\frac{1}{2}$, $\Real\kappa<\frac{1}{2}$, 
$-1<t<1$, $-1<x<1$.
\end{abstract}
%
%\newpage
%
\section{Introduction}
\label{I}
\setcounter{equation}{0}
Recently, Cohl \cite{Cohl11} has derived the integral formula
\begin{eqnarray}
&& \hspace*{-2em}
\frac{n+\lambda}{\lambda}\int_{-1}^{1}\mathrm{d}t\:
\frac{(1-t^{2})^{\lambda-1/2}C_{n}^{\lambda}(t)}{(z-t)^{\kappa+1/2}}
=\mathrm{e}^{\mathrm{i}\pi(\lambda-\kappa)}
\frac{\sqrt{\pi}\,(n+\lambda)\Gamma(n+2\lambda)}
{2^{\lambda-3/2}n!\Gamma(\lambda+1)\Gamma(\kappa+\frac{1}{2})}
(z^{2}-1)^{(\lambda-\kappa)/2}
\mathfrak{Q}_{n+\lambda-1/2}^{\kappa-\lambda}(z)
\nonumber \\
&& \hspace*{20em} (\textrm{$\Real\lambda>{\textstyle-\frac{1}{2}}$,
$\kappa\in\mathbb{C}$, $z\in\mathbb{C}\setminus(-\infty,1]$}),
\label{1.1}
\end{eqnarray}
where $C_{n}^{\lambda}(t)$ is the Gegenbauer polynomial, while
$\mathfrak{Q}_{\nu}^{\mu}(z)$ is the associated Legendre function of
the second kind. The integral (\ref{1.1}) generalizes the Gormley's
result \cite{Gorm34}
\begin{eqnarray}
&& \hspace*{-2em}
\frac{n+\lambda}{\lambda}\int_{-1}^{1}\mathrm{d}t\:
\frac{(1-t^{2})^{\lambda-1/2}C_{n}^{\lambda}(t)}{z-t}
=\mathrm{e}^{\mathrm{i}\pi(\lambda-1/2)}
\frac{\sqrt{\pi}\,(n+\lambda)\Gamma(n+2\lambda)}
{2^{\lambda-3/2}n!\Gamma(\lambda+1)}(z^{2}-1)^{(\lambda-1/2)/2}
\mathfrak{Q}_{n+\lambda-1/2}^{1/2-\lambda}(z)
\nonumber \\
&& \hspace*{20em} (\textrm{$\Real\lambda>{\textstyle-\frac{1}{2}}$,
$z\in\mathbb{C}\setminus(-\infty,1]$}),
\label{1.2}
\end{eqnarray}
which, in turn, is an extension of the celebrated Neumann's integral
formula \cite{Neum48}
\begin{equation}
\int_{-1}^{1}\mathrm{d}t\:\frac{P_{n}(t)}{z-t}=2\mathfrak{Q}_{n}(z)
\qquad (z\in\mathbb{C}\setminus[-1,1]),
\label{1.3}
\end{equation}
where $P_{n}(t)$ is the Legendre polynomial. The factor
$(n+\lambda)/\lambda$ appearing in front of the integrals in Eqs.\
(\ref{1.1}) and (\ref{1.2}) (and also at some other places in the
text) is only seemingly awkward and has been introduced to avoid the
difficulty one otherwise encounters when $\lambda\to0$.

From Eq.\ (\ref{1.1}) and from the closure relation for the
Gegenbauer polynomials, which is
\begin{eqnarray}
&& \frac{2^{2\lambda-1}\Gamma^{2}(\lambda)}{\pi}
\sum_{n=0}^{\infty}\frac{n!(n+\lambda)}{\Gamma(n+2\lambda)}
C_{n}^{\lambda}(t)C_{n}^{\lambda}(t')
=\frac{\delta(t-t')}{(1-t^{2})^{(\lambda-1/2)/2}
(1-t^{\prime\,2})^{(\lambda-1/2)/2}}
\nonumber \\
&& \hspace*{20em} (\textrm{$\Real\lambda>{\textstyle-\frac{1}{2}}$, 
$-1<t,t'<1$})
\label{1.4}
\end{eqnarray}
(here $\delta(t-t')$ is the Dirac delta function), one may deduce the
summation formula \cite{Cohl11}
\begin{eqnarray}
&& \sum_{n=0}^{\infty}\frac{n+\lambda}{\lambda}
\mathfrak{Q}_{n+\lambda-1/2}^{\kappa-\lambda}(z)C_{n}^{\lambda}(t)
=\mathrm{e}^{\mathrm{i}\pi(\kappa-\lambda)}
\frac{\sqrt{\pi}\,\Gamma(\kappa+\frac{1}{2})}
{2^{\lambda+1/2}\Gamma(\lambda+1)}
\frac{(z^{2}-1)^{(\kappa-\lambda)/2}}{(z-t)^{\kappa+1/2}}
\nonumber \\
&& \hspace*{10em} (\textrm{$\Real\lambda>{\textstyle-\frac{1}{2}}$, 
$\kappa\in\mathbb{C}$, $-1<t<1$, 
$z\in\mathbb{C}\setminus(-\infty,1]$}). 
\label{1.5} 
\end{eqnarray}
The particular case of this relation with $\kappa=\lambda$ has been
known before \cite[p.\ 183]{Magn66}. For $\kappa=\lambda=1/2$, Eq.\
(\ref{1.5}) reduces to the Heine's identity
\begin{equation}
\sum_{n=0}^{\infty}(2n+1)\mathfrak{Q}_{n}(z)P_{n}(t)=\frac{1}{z-t},
\qquad
(\textrm{$-1<t<1$, $z\in\mathbb{C}\setminus(-\infty,1]$}).
\label{1.6}
\end{equation}

In Sec.\ \ref{II} of this communication, we show that one may use the
relation in Eq.\ (\ref{1.1}) to evaluate some further definite
integrals involving the Gegenbauer polynomials. Furthermore, in Sec.\
\ref{III}, we exploit the identity (\ref{1.5}) to determine
closed-form representations of some series involving the Gegenbauer
polynomials and the associated Legendre functions on the cut. While
particular cases of the relations we arrive at in the present work,
corresponding to specific choices of the parameters $\lambda$ and
$\kappa$, may be found in the literature on special functions, we are
not aware of any appearance of these relations in their most general
forms derived below.

Throughout the paper, it is understood that
\begin{equation}
(z^{2}-1)^{\alpha}\equiv(z-1)^{\alpha}(z+1)^{\alpha}
\qquad (|\arg(z\pm1)|<\pi).
\label{1.7}
\end{equation}
%
%\newpage
%
\section{Evaluation of some definite integrals involving the
Gegenbauer polynomials}
\label{II}
\setcounter{equation}{0}
Let us assume that $\Real\kappa<\frac{1}{2}$. We proceed to the
investigation of the limit of the formula in Eq.\ (\ref{1.1}) as
$z\to x\pm\mathrm{i}0$, with $-1<x<1$. Exploiting the identities
\begin{equation}
x\pm\mathrm{i}0+1=1+x,
\qquad
x\pm\mathrm{i}0-1=\mathrm{e}^{\pm\mathrm{i}\pi}(1-x)
\qquad (-1<x<1)
\label{2.1}
\end{equation}
and
\begin{equation}
\mathrm{e}^{-\mathrm{i}\pi\mu}
\mathfrak{Q}_{\nu}^{\mu}(x\pm\mathrm{i}0)
=\mathrm{e}^{\pm\mathrm{i}\pi\mu/2}\left[Q_{\nu}^{\mu}(x)
\mp\frac{\mathrm{i}\pi}{2}P_{\nu}^{\mu}(x)\right]
\qquad (-1<x<1),
\label{2.2}
\end{equation}
where $P_{\nu}^{\mu}(x)$ and $Q_{\nu}^{\mu}(x)$ are the associated
Legendre functions \emph{on the cut\/}, we obtain
\begin{eqnarray}
\frac{n+\lambda}{\lambda}\int_{-1}^{1}\mathrm{d}t\:
\frac{(1-t^{2})^{\lambda-1/2}C_{n}^{\lambda}(t)}
{(x\pm\mathrm{i}0-t)^{\kappa+1/2}}
&=& \frac{\sqrt{\pi}\,(n+\lambda)\Gamma(n+2\lambda)}
{2^{\lambda-3/2}n!\Gamma(\lambda+1)\Gamma(\kappa+\frac{1}{2})}
\nonumber \\
&& \times\,(1-x^{2})^{(\lambda-\kappa)/2}
\left[Q_{n+\lambda-1/2}^{\kappa-\lambda}(x)
\mp\frac{\mathrm{i}\pi}{2}P_{n+\lambda-1/2}^{\kappa-\lambda}(x)\right]
\nonumber \\
&& \hspace*{5em} (\textrm{$\Real\lambda>{\textstyle-\frac{1}{2}}$, 
$\Real\kappa<{\textstyle\frac{1}{2}}$, $-1<x<1$}).
\label{2.3}
\end{eqnarray}
(The reason for imposing the above constraint on $\kappa$ is now
clear: we have had to exclude a non-integrable singularity at $t=x$
that otherwise occurs.) Since it holds that
\begin{equation}
x\pm\mathrm{i}0-t
=\left\{
\begin{array}{ll}
x-t & (t<x) \\
\mathrm{e}^{\pm\mathrm{i}\pi}(t-x) & (t>x),
\end{array}
\right.
\label{2.4}
\end{equation}
we may split the integrals appearing on the left-hand side of Eq.\
(\ref{2.3}) as follows:
\begin{eqnarray}
\frac{n+\lambda}{\lambda}\int_{-1}^{1}\mathrm{d}t\:
\frac{(1-t^{2})^{\lambda-1/2}C_{n}^{\lambda}(t)}
{(x\pm\mathrm{i}0-t)^{\kappa+1/2}}
&=& \frac{n+\lambda}{\lambda}\int_{-1}^{x}\mathrm{d}t\:
\frac{(1-t^{2})^{\lambda-1/2}C_{n}^{\lambda}(t)}{(x-t)^{\kappa+1/2}}
\nonumber \\
&& +\,\mathrm{e}^{\mp\mathrm{i}\pi(\kappa+1/2)}
\frac{n+\lambda}{\lambda}\int_{x}^{1}\mathrm{d}t\:
\frac{(1-t^{2})^{\lambda-1/2}C_{n}^{\lambda}(t)}{(t-x)^{\kappa+1/2}}
\nonumber \\
&& (\textrm{$\Real\lambda>{\textstyle-\frac{1}{2}}$,
$\Real\kappa<{\textstyle\frac{1}{2}}$, $-1<x<1$}).
\label{2.5}
\end{eqnarray}
If we insert Eq.\ (\ref{2.5}) into Eq.\ (\ref{2.3}), and then equate
separately terms corresponding to the choices of the upper and the
lower signs on both sides of the resulting equation, we arrive at an
inhomogeneous algebraic system for the two integrals standing on the
right-hand side of Eq.\ (\ref{2.5}). Solving this system, and using
the well-known identity
\begin{equation}
\Gamma(\zeta)\Gamma(1-\zeta)=\frac{\pi}{\sin(\pi\zeta)},
\label{2.6}
\end{equation}
we obtain
\begin{eqnarray}
&& \hspace*{-2em}
\frac{n+\lambda}{\lambda}\int_{x}^{1}\mathrm{d}t\:
\frac{(1-t^{2})^{\lambda-1/2}C_{n}^{\lambda}(t)}{(t-x)^{\kappa+1/2}}
=\frac{\sqrt{\pi}\,(n+\lambda)\Gamma(n+2\lambda)
\Gamma(\frac{1}{2}-\kappa)}{2^{\lambda-1/2}n!\Gamma(\lambda+1)}
(1-x^{2})^{(\lambda-\kappa)/2}P_{n+\lambda-1/2}^{\kappa-\lambda}(x)
\nonumber \\
&& \hspace*{20em} (\textrm{$\Real\lambda>{\textstyle-\frac{1}{2}}$,
$\Real\kappa<{\textstyle\frac{1}{2}}$, $-1<x<1$})
\label{2.7}
\end{eqnarray}
and
\begin{eqnarray}
&& \hspace*{-2em}
\frac{n+\lambda}{\lambda}\int_{-1}^{x}\mathrm{d}t\:
\frac{(1-t^{2})^{\lambda-1/2}C_{n}^{\lambda}(t)}{(x-t)^{\kappa+1/2}}
\nonumber \\
&& =\,\frac{\sqrt{\pi}\,(n+\lambda)\Gamma(n+2\lambda)}
{2^{\lambda-3/2}n!\Gamma(\lambda+1)\Gamma(\kappa+\frac{1}{2})}
(1-x^{2})^{(\lambda-\kappa)/2}
\left\{Q_{n+\lambda-1/2}^{\kappa-\lambda}(x)
-\frac{\pi}{2}P_{n+\lambda-1/2}^{\kappa-\lambda}(x)
\cot[\pi(\kappa+{\textstyle\frac{1}{2}})]\right\}
\nonumber \\
&& \hspace*{18em} (\textrm{$\Real\lambda>{\textstyle-\frac{1}{2}}$,
$\Real\kappa<{\textstyle\frac{1}{2}}$, $-1<x<1$}),
\label{2.8}
\end{eqnarray}
respectively. The right-hand side of the latter equation may be
simplified considerably after one exploits the known
relation\footnote{~It is stated in Refs.\ \cite[p.\ 170]{Magn66} and
\cite[p.\ 144]{Erde53a} that the domain of validity of the relation
displayed in our Eq.\ (\ref{2.9}), and also of the counterpart
expression for $Q_{\nu}^{\mu}(-x)$ in terms of $P_{\nu}^{\mu}(x)$ and
$Q_{\nu}^{\mu}(x)$, is $0<x<1$. However, it is not difficult to show
that if both relations hold on that interval, they must be valid for
$-1<x\leqslant0$ as well.}
\begin{equation}
P_{\nu}^{\mu}(-x)=P_{\nu}^{\mu}(x)\cos[\pi(\nu+\mu)]
-\frac{2}{\pi}Q_{\nu}^{\mu}(x)\sin[\pi(\nu+\mu)]
\qquad (-1<x<1)
\label{2.9}
\end{equation}
and the identity (\ref{2.6}). This yields
\begin{eqnarray}
&& \frac{n+\lambda}{\lambda}\int_{-1}^{x}\mathrm{d}t\:
\frac{(1-t^{2})^{\lambda-1/2}C_{n}^{\lambda}(t)}{(x-t)^{\kappa+1/2}}
\nonumber \\
&& \qquad =\,(-)^{n}\frac{\sqrt{\pi}\,(n+\lambda)\Gamma(n+2\lambda)
\Gamma(\frac{1}{2}-\kappa)}{2^{\lambda-1/2}n!\Gamma(\lambda+1)}
(1-x^{2})^{(\lambda-\kappa)/2}P_{n+\lambda-1/2}^{\kappa-\lambda}(-x)
\nonumber \\
&& \hspace*{15em} (\textrm{$\Real\lambda>{\textstyle-\frac{1}{2}}$,
$\Real\kappa<{\textstyle\frac{1}{2}}$, $-1<x<1$}).
\label{2.10}
\end{eqnarray}
Equation (\ref{2.10}) may be also derived from Eq.\ (\ref{2.7}) after
in the latter one makes the simultaneous replacements $t\to-t$ and
$x\to-x$, and subsequently exploits the property
$C_{n}^{\lambda}(-t)=(-)^{n}C_{n}^{\lambda}(t)$.

Equations (\ref{2.7}) and (\ref{2.10}) constitute the result of this
section. Their particular cases, corresponding to making the choices
$\lambda=\frac{1}{2}$ and $\kappa=0$, may be found, albeit with the
right-hand sides written in less compact forms, in Refs.\ \cite[p.\
261]{Magn66} and \cite[p.\ 187]{Erde53b}.
%
%\newpage
%
\section{Evaluation of closed forms of some series involving the
Gegenbauer polynomials and the associated Legendre functions on the
cut}
\label{III}
\setcounter{equation}{0}
Next, we shall draw inferences from the expansion (\ref{1.5}),
approaching the limit $z\to x\pm\mathrm{i}0$, with $-1<x<1$. As in
the preceding section, we assume that $\Real\kappa<\frac{1}{2}$
(recall that the formula in Eq.\ (\ref{1.5}) is a consequence of Eq.\
(\ref{1.1}) which in the limit discussed here loses its sense unless
the above restriction is imposed on $\kappa$). With the use of Eqs.\
(\ref{2.1}), (\ref{2.2}) and (\ref{2.4}), we arrive at
\begin{eqnarray}
&& \sum_{n=0}^{\infty}\frac{n+\lambda}{\lambda}
\left[Q_{n+\lambda-1/2}^{\kappa-\lambda}(x)\mp\frac{\mathrm{i}\pi}{2}
P_{n+\lambda-1/2}^{\kappa-\lambda}(x)\right]C_{n}^{\lambda}(t)
\nonumber \\
&& \qquad=\,\frac{\sqrt{\pi}\,\Gamma(\kappa+\frac{1}{2})}
{2^{\lambda+1/2}\Gamma(\lambda+1)}(1-x^{2})^{(\kappa-\lambda)/2}
\times\left\{
\begin{array}{lc}
(x-t)^{-\kappa-1/2} & (-1<t<x<1) \\
\mathrm{e}^{\mp\mathrm{i}\pi(\kappa+1/2)}
(t-x)^{-\kappa-1/2} & (-1<x<t<1)
\end{array}
\right.
\nonumber \\
&& \hspace*{25em} (\textrm{$\Real\lambda>{\textstyle-\frac{1}{2}}$,
$\Real\kappa<{\textstyle\frac{1}{2}}$}).
\label{3.1}
\end{eqnarray}
Hence, subtracting or adding the two relations embodied in Eq.\
(\ref{3.1}), we deduce that
\begin{eqnarray}
&& \sum_{n=0}^{\infty}\frac{n+\lambda}{\lambda}
P_{n+\lambda-1/2}^{\kappa-\lambda}(x)C_{n}^{\lambda}(t)
\nonumber \\
&& \qquad =\,\frac{\sqrt{\pi}}
{2^{\lambda-1/2}\Gamma(\lambda+1)\Gamma(\frac{1}{2}-\kappa)}
(1-x^{2})^{(\kappa-\lambda)/2}
\times\left\{
\begin{array}{lc}
0 & (-1<t<x<1) \\
(t-x)^{-\kappa-1/2} & (-1<x<t<1)
\end{array}
\right.
\nonumber \\
&& \hspace*{25em} (\textrm{$\Real\lambda>{\textstyle-\frac{1}{2}}$,
$\Real\kappa<{\textstyle\frac{1}{2}}$})
\label{3.2}
\end{eqnarray}
and
\begin{eqnarray}
&& \sum_{n=0}^{\infty}\frac{n+\lambda}{\lambda}
Q_{n+\lambda-1/2}^{\kappa-\lambda}(x)C_{n}^{\lambda}(t)
\nonumber \\
&& \qquad =\,\frac{\sqrt{\pi}\,\Gamma(\kappa+\frac{1}{2})}
{2^{\lambda+1/2}\Gamma(\lambda+1)}(1-x^{2})^{(\kappa-\lambda)/2}
\times\left\{
\begin{array}{lc}
(x-t)^{-\kappa-1/2} & (-1<t<x<1) \\
(t-x)^{-\kappa-1/2}\cos[\pi(\kappa+{\textstyle\frac{1}{2}})] 
& (-1<x<t<1)
\end{array}
\right.
\nonumber \\
&& \hspace*{25em} (\textrm{$\Real\lambda>{\textstyle-\frac{1}{2}}$,
$\Real\kappa<{\textstyle\frac{1}{2}}$}).
\label{3.3}
\end{eqnarray}
If in Eqs.\ (\ref{3.2}) and (\ref{3.3}) we make the simultaneous
replacements $t\to-t$ and $x\to-x$, we obtain two further summation
formulas:
\begin{eqnarray}
&& \sum_{n=0}^{\infty}(-)^{n}\frac{n+\lambda}{\lambda}
P_{n+\lambda-1/2}^{\kappa-\lambda}(-x)C_{n}^{\lambda}(t)
\nonumber \\
&& \qquad =\,\frac{\sqrt{\pi}}
{2^{\lambda-1/2}\Gamma(\lambda+1)\Gamma(\frac{1}{2}-\kappa)}
(1-x^{2})^{(\kappa-\lambda)/2}
\times\left\{
\begin{array}{lc}
(x-t)^{-\kappa-1/2} & (-1<t<x<1) \\
0 & (-1<x<t<1)
\end{array}
\right.
\nonumber \\
&& \hspace*{25em} (\textrm{$\Real\lambda>{\textstyle-\frac{1}{2}}$,
$\Real\kappa<{\textstyle\frac{1}{2}}$}),
\label{3.4}
\end{eqnarray}
\begin{eqnarray}
&& \sum_{n=0}^{\infty}(-)^{n}\frac{n+\lambda}{\lambda}
Q_{n+\lambda-1/2}^{\kappa-\lambda}(-x)C_{n}^{\lambda}(t)
\nonumber \\
&& \qquad =\,\frac{\sqrt{\pi}\,\Gamma(\kappa+\frac{1}{2})}
{2^{\lambda+1/2}\Gamma(\lambda+1)}(1-x^{2})^{(\kappa-\lambda)/2}
\times\left\{
\begin{array}{lc}
(x-t)^{-\kappa-1/2}\cos[\pi(\kappa+{\textstyle\frac{1}{2}})] 
& (-1<t<x<1) \\
(t-x)^{-\kappa-1/2} & (-1<x<t<1)
\end{array}
\right.
\nonumber \\
&& \hspace*{25em} (\textrm{$\Real\lambda>{\textstyle-\frac{1}{2}}$,
$\Real\kappa<{\textstyle\frac{1}{2}}$}).
\label{3.5}
\end{eqnarray}
Some particular cases of the expansions (\ref{3.1})--(\ref{3.5}),
corresponding to specific choices of $\kappa$ and/or $\lambda$, may
be found\footnote{~In the second series in Sec.\ 4.5.4 in Ref.\
\cite[p.\ 182]{Magn66}, $P_{m-1/2}^{\mu}(\cos\vartheta)$ should be
replaced by $P_{m-1/2}^{\mu}(-\cos\vartheta)$. In the first formula
in Ref.\ \cite[p.\ 183]{Magn66}, the constraint $x<\cos\varphi$
should be replaced by $-1<\cos\varphi<x<1$. Moreover, it follows from
our Eq.\ (\ref{3.4}) that the latter formula is valid at least for
$\Real\nu>-\frac{1}{2}$.} in Refs.\ \cite[pp.\ 182--183]{Magn66},
\cite[p.\ 166]{Erde53a} and \cite[pp.\ 341--342]{Prud03}.
\section*{Acknowledgments}
I thank Dr.\ Howard S.\ Cohl for drawing my attention to his work
\cite{Cohl11} and for a stimulating correspondence on the Legendre
functions and related subjects.
\end{document}